\documentclass{article}
\usepackage{amssymb}
\usepackage{amsbsy}
\begin{document}
\centerline{\Large \bf } \vskip 6pt

\begin{center}{\Large \bf Orlicz log-Aleksandrov-Fenchel inequality}\end{center}

\vskip 6pt\begin{center} \centerline{Chang-Jian
Zhao} \centerline{\it Department of
Mathematics, China Jiliang University, Hangzhou 310018, P. R.
China}\centerline{\it Email: chjzhao@aliyun.com}
\end{center}

\begin{center}
\begin{minipage}{12cm}
{\bf Abstract} In this paper, we establish an Orlicz log-Aleksandrov-Fenchel inequality by introducing new concepts of mixed volume measure and Orlicz multiple mixed volume measure, and using the Orlicz-Aleksandrov-Fenchel inequality. The Orlicz log-Aleksandrov-Fenchel inequality in special cases yield the classical Aleksandrov-Fenchel inequality and Orlicz log-Minkowski type inequality, respectively.

{\bf Keywords} mixed volumes, $L_{p}$-mixed volume, Minkowski inequality, logarithmic Minkowski inequality, Aleksandrov-Fenchel inequality, Orlicz-Aleksandrov-Fenchel inequality.

{\bf 2010 Mathematics Subject Classification} 46E30 52A39
\end{minipage}
\end{center}
\vskip 20pt

\vskip 10pt \noindent{\large \bf 1 ~Introducation}\vskip 10pt

In 2012, B\"{o}r\"{o}czky, Lutwak, Yang, and Zhang [1] conjecture
a {\it logarithmic Minkowski inequality} for origin-symmetric convex
bodies, the logarithmic Minkowski inequality was stated the following.

{\bf The conjectured logarithmic Minkowski inequality}~ {\it If $K$ and $L$ are convex bodies in ${\Bbb R}^n$ which are symmetric with respect to the origin, then
$$\int_{S^{n-1}}\ln\left(\frac{h_{K}}{h_{L}}\right)d\bar{V}_{L}
\geq\frac{1}{n}\ln\left(\frac{V(K)}{V(L)}\right),\eqno(1.1)$$
where $dv_{L}=\frac{1}{n}h_{L}dS_{L}$ is the cone-volume measure of $L$, and
$d\bar{V}_{L}=\frac{1}{V(L)}dv_{L}$ is its normalization,
and $S_{L}=S(L,\cdot)$ is the mixed surface area measure of} $L$.

The functions are the support functions. If $K$ is
a nonempty closed (not necessarily bounded) convex set in ${\Bbb
R}^{n}$, then
$$h_{K}=\max\{x\cdot y: y\in K\},$$
for $x\in {\Bbb R}^{n},$ defines the support function $h_{K}$ of
$K$. A nonempty closed convex set is uniquely determined by its
support function.

In 2016, Stancu [2] proved a modified logarithmic Minkowski inequality for non-symmetric convex bodies not symmetric with respect to the origin. The logarithmic Minkowski inequality was given in the following theorem.

{\bf The logarithmic Minkowski inequality}~ {\it If $K$ and $L$ are convex bodies in ${\Bbb R}^{n}$ that containing the origin in their interior, then
$$\int_{S^{n-1}}\ln\left(\frac{h_{K}}{h_{L}}\right)d\bar{V}_{1}
\geq\frac{1}{n}\ln\left(\frac{V(K)}{V(L)}\right).\eqno(1.2)$$
with equality if and only if $K$ and $L$ are homothetic, where $dv_{1}$ is the mixed volume measure
$dv_{1}=\frac{1}{n}h_{K}dS_{L},$ and $d\bar{v}_{1}=\frac{1}{V_{1}(L,K)}dv_{1}$
is its normalization, and $V_{1}(L,K)$ denotes the usual mixed volume of $L$ and $K$, defined by} ({\it see} []) $$V_{1}(L,K)=\frac{1}{n}\int_{S^{n-1}}h_{K}dS_{L}.$$

Recently, this logarithmic Minkowski inequality and the conjectured logarithmic Minkowski inequality
have attracted extensive attention and research. The recent research on the logarithmic Minkowski and its dual type inequalities can be found in the references [3-16].

Associated with the convex bodies $K_{1},\cdots,K_{n-1}$ in ${\Bbb
R}^{n}$ is a unique positive Borel measure on $S^{n-1}$,
$S(K_{1},\ldots,K_{n-1};\cdot)$, call it the mixed area measure of
$K_{1},\ldots,K_{n-1}$, with the property that any convex body
$K_{n}$ one has the integral representation (see e.g. [17], p.
354).
$$V(K_{1},\ldots,K_{n})=\frac{1}{n}\int_{S^{n-1}}h_{K_{n}}dS(K_{1},\ldots,K_{n-1};u).\eqno(1.3)$$
The integration is with respect to the mixed area measure
$S(K_{1},\ldots,K_{n-1};\cdot)$ on $S^{n-1}$. The mixed area
measure $S(K_{1},\ldots,K_{n-1};\cdot)$ is symmetric in its (first
$n-1$) arguments. When $K_{1}=\cdots=K_{n-i-1}=K$ and
$K_{n-i}=\cdots=K_{n-1}=B$, the mixed area measure
$S(K,\ldots,K,B,\ldots,B;\cdot)$ with $i$ copies of $B$ and
$(n-i-1)$ copies of $K$, will be written as $S_{i}(K,\cdot)$. For
$K_{1}=\cdots=K_{n-1}=K$, $S(K_{1},\ldots,K_{n-1};\cdot)$ reduces
to the surface area measure $S(K;\cdot)$.

It is well known that in Brunn-Minkowski theory, Minkowski inequality and Aleksandrov-Fenchel inequality appear at the same time, and the latter is a generalization of the former. So a natural question is raised: is there a logarithmic Aleksandrov-Fenchel inequality relative to the logarithmic Minkowski inequality? The main purpose of this article is to answer the above questions perfectly and obtain a logarithmic Aleksandrov-Fenchel inequality  by introducing two new concepts of mixed volume measure and $L_{p}$-multiple mixed volume measure, and using the $L_{p}$-Aleksandrov-Fenchel inequality for the $L_{p}$-multiple mixed volume. The logarithmic Aleksandrov-Fenchel inequality in special cases yield  the classical Aleksandrov-Fenchel inequality, and four recent logarithmic Minkowski type inequalities, which are logarithmic Minkowski inequality for mixed volumes, $L_{p}$-mixed volumes, quermassintegrals and $p$-mixed quermassintergrals, respectively. Our main result is given in the following theorem.

{\bf Orlicz log-Aleksandrov-Fenchel inequality}~ {\it Let $L_{1},\ldots,L_{n},K_{n}$ be convex bodies in ${\Bbb R}^{n}$ that containing the origin in their interior. If $\varphi:[0,\infty)\rightarrow(0,\infty)$ is a convex and
increasing function such that $\varphi(0)=0$ and $\varphi(1)=1$, and $1\leq r\leq n$, then}
$$\int_{S^{n-1}}\ln\left(\varphi\left(\frac{h_{K_{n}}}{{h_{L_{n}}}}\right)\right)
d\bar{V}_{\varphi}(L_{1},\ldots,L_{n-1},K_{n},L_{n})\geq\ln\left(\varphi\left(\frac{\prod_{i=1}^{r}V(L_{i},\ldots,L_{i},L_{r+1},\ldots,L_{n-1},K_{n})^{1/r}}
{V(L_{1},\ldots,L_{n})}\right)\right).\eqno(1.4)$$

Here $d\bar{V}_{\varphi}(L_{1},\ldots,L_{n-1},K_{n},L_{n})$ denotes a new probability measure call it Orlicz multiple mixed volume probability measure of convex bodies $L_{1},\ldots,L_{n},K_{n}$, defined by
$$d\bar{V}_{\varphi}(L_{1},\cdots,L_{n-1},K_{n},L_{n})=\frac{dv_{\varphi}(L_{1},\cdots,L_{n-1},K_{n},L_{n})}{V_{\varphi}(L_{1},
\cdots,L_{n-1},K_{n},L_{n})},\eqno(1.5)$$
and $dv_{\varphi}(L_{1},\cdots,L_{n-1},K_{n},L_{n})$ denotes the Orlicz multiple mixed volume measure of $L_{1},\ldots,L_{n},K_{n}$, defined by
$$dv_{\varphi}(L_{1},\cdots,L_{n-1},K_{n},L_{n})=\frac{1}{n}\cdot\varphi\left(\frac{h_{K_{n}}}{h_{L_{n}}}
\right)h_{L_{n}}dS(L_{1},\ldots,L_{n-1};\cdot).\eqno(1.6)$$
Moreover, $V_{\varphi}(L_{1},\ldots,L_{n-1},K_{n},L_{n})$ is the Orlicz multiple mixed volume of convex bodies $L_{1},\ldots,L_{n},K_{n}$ defined by ([18])
$$V_{\varphi}(L_{1},\cdots,L_{n-1},K_{n},L_{n})=\frac{1}{n}\int_{S^{n-1}}
\varphi\left(\frac{h_{L_{n}}}{h_{K_{n}}}\right)h_{K_{n}}dS(L_{1},\ldots,L_{n-1};u).\eqno(1.7)$$

{\bf Remark}~ When $K_{n}=L_{n}$, inequality (1.4) becomes the classical Aleksandrov-Fenchel inequality of convex bodies $L_{1},\ldots,L_{n}$ as follows: If
$L_{1},\cdots,L_{n}$ are convex bodies that containing the origin and
$1\leq r\leq n$, then
$$V(L_{1},\cdots,
L_{n})\geq\prod_{i=1}^{r}V(L_{i}\ldots,L_{i},L_{r+1},\ldots,L_{n})^{1/r},\eqno(1.8)$$ (see e.g. [19, p. 401]).

On the other hand, when $\varphi(x)=x^{p}$, $p=1$, $r=n-1$, $L_{1}=\ldots=L_{n-1}=L$, $L_{n}=L$ and $K_{n}=K$, then (1.4) becomes Stancu's logarithmic Minkowski inequality (1.2) established in [2].

\vskip 10pt \noindent{\large \bf 2 ~Notations and
preliminaries}\vskip 10pt

The setting for this paper is $n$-dimensional Euclidean space
${\Bbb R}^{n}$. Let ${\cal K}^{n}$ be the class of nonempty
compact convex subsets of ${\Bbb R}^{n}$, let ${\cal K}^{n}_{o}$
be the class of members of ${\cal K}^{n}$ containing the origin. A set $K\in {\cal
K}^{n}$ is called a convex body if its interior is nonempty. We
reserve the letter $u\in S^{n-1}$ for unit vectors, and the letter
$B$ for the unit ball centered at the origin. The surface of $B$
is $S^{n-1}$. For a compact set $K$, we write $V(K)$ for the
($n$-dimensional) Lebesgue measure of $K$ and call this the volume
of $K$. Let $d$ denote the Hausdorff metric on ${\cal K}^{n}$, i.e., for
$K, L\in {\cal K}^{n},$ $$d(K,L)=|h_{K}-h_{L}|_{\infty},$$ where
$|\cdot|_{\infty}$ denotes the sup-norm on the space of continuous
functions $C(S^{n-1})$.

\vskip 10pt{\it ~2.1 Mixed quermassintegrals}\vskip 10pt

If $K_{i}\in {\cal K}^{n}$ $(i=1,2,\ldots,r)$ and
     $\lambda_{i}$ $(i=1,2,\ldots,r)$ are nonnegative real numbers, then
     of fundamental importance is the fact that the volume of
     $\sum_{i=1}^{r}\lambda_{i}K_{i}$ is a homogeneous polynomial
     in $\lambda_{i}$ given by (see e.g. [15])
     $$V(\lambda_{1}K_{1}+\cdots+\lambda_{n}K_{n})=\sum_{i_{1},\ldots,i_{n}}
     \lambda_{i_{1}}\ldots\lambda_{i_{n}}V_{i_{1}\ldots
     i_{n}},\eqno(2.1)$$ where the sum is taken over all
     $n$-tuples $(i_{1},\ldots,i_n)$ of positive integers not
     exceeding $r$. The coefficient $V_{i_{1}\ldots i_n}$ depends
     only on the bodies $K_{i_{1}},\ldots,K_{i_{n}}$ and is
     uniquely determined by (2.1), it is called the mixed volume
     of $K_{i},\ldots,K_{i_{n}}$, and is written as
     $V(K_{i_{1}},\ldots,K_{i_{n}}).$
     Let $K_{1}=\ldots=K_{n-i}=K$ and $K_{n-i+1}=\ldots=K_{n}=L$, then
the mixed volume $V(K_{1},\ldots, K_{n})$ is written as $V_{i}(K,L)$. If $K_{1}=\cdots=K_{n-i}=K,$
$K_{n-i+1}=\cdots=K_{n}=B$ The mixed volumes $V_{i}(K,B)$ is written as $W_{i}(K)$ and call as quermassintegrals (or
$i$th mixed quermassintegrals) of $K$. We write $W_{i}(K,L)$ for
the mixed volume $V(\underbrace{K,\ldots,K}_{n-i-1},\underbrace{B,\ldots,B}_{i},L)$ and call as mixed
quermassintegrals of convex bodies $K$ and $L$. Aleksandrov [20] and Fenchel and Jessen [21]
(also see Busemann [22] and Schneider [23]) have shown that for
$K\in {\cal K}^{n}_{o}$, and $i=0,1,\ldots,n-1,$ there exists a
regular Borel measure $S_{i}(K,\cdot)$ on $S^{n-1}$, such that the
mixed quermassintegrals $W_{i}(K,L)$ has the following
representation:
$$W_{i}(K,L)=\frac{1}{n-i}\lim_{\varepsilon\rightarrow 0^{+}}\frac{W_{i}(K+\varepsilon L)-W_{i}(K)}{\varepsilon}=\frac{1}{n}
\int_{S^{n-1}}h_{L}dS_{i}(K,u).\eqno(2.2)$$ Associated with
$K_{1},\ldots, K_{n}\in {\cal K}^{n}$ is a Borel measure
$S(K_{1},\ldots,K_{n-1},\cdot)$ on $S^{n-1}$, called the mixed
surface area measure of $K_{1},\ldots, K_{n-1}$, which has the
property that for each $K\in {\cal K}^{n}$ (see e.g. [24], p.353),
$$V(K_{1},\ldots,K_{n-1},K)=\frac{1}{n}\int_{S^{n-1}}h_{K}dS(K_{1},\ldots,K_{n-1},u).\eqno(2.3)$$
In fact, the measure $S(K_{1},\ldots,K_{n-1},\cdot)$ can be
defined by the propter that (2.3) holds for all $K\in {\cal
K}^{n}.$ Let $K_{1}=\ldots=K_{n-i-1}=K$ and
$K_{n-i}=\ldots=K_{n-1}=L$, then the mixed surface area measure
$S(K_{1},\ldots, K_{n-1},\cdot)$ is written as $S_{i}(K,L;\cdot)$. When $L=B$, $S_{i}(K,L;\cdot)$ is written as
$S_{i}(K,\cdot)$ and called as $i$th mixed surface area measure. A
fundamental inequality for mixed quermassintegrals stats that: For
$K,L\in{\cal K}^{n}$ and $0\leq i<n-1$,
$$W_{i}(K,L)^{n-i}\geq W_{i}(K)^{n-i-1}W_{i}(L),\eqno(2.4)$$
with equality if and only if $K$ and $L$ are homothetic and
$L=\{o\}$.

\vskip 8pt {\it 2.2 ~Mixed $p$-quermassintegrals}\vskip 10pt

Mixed quermassintegrals are, of course, the first variation of the
ordinary quermassintegrals, with respect to Minkowski addition.
The mixed quermassintegrals $W_{p,0}(K,L),
W_{p,1}(K,L),\ldots,W_{p,n-1}(K,L)$, as the first variation of the
ordinary quermassintegrals, with respect to Firey addition: For
$K,L\in{\cal K}_{o}^{n}$, and real $p\geq 1$, defined by (see
e.g. [25])
$$W_{p,i}(K,L)=\frac{p}{n-i}\lim_{\varepsilon\rightarrow _{0^{+}}}\frac{W_{i}(K+_{p}\varepsilon\cdot L)-W_{i}(K)}{\varepsilon}.\eqno(2.5)$$
The mixed $p$-quermassintegrals $W_{p,i}(K,L)$, for all $K,L\in
{\cal K}^{n}_{o}$, has the following integral representation:
$$W_{p,i}(K,L)=\frac{1}{n}\int_{S^{n-1}}h_{L}^{p}dS_{p,i}(K,u),\eqno(2.6)$$
where $S_{p,i}(K,\cdot)$ denotes the Boel measure on $S^{n-1}$.
The measure $S_{p,i}(K,\cdot)$ is absolutely continuous with
respect to $S_{i}(K,\cdot)$, and has Radon-Nikodym derivative
$$\frac{dS_{p,i}(K,\cdot)}{dS_{i}(K,\cdot)}=h_{K}^{1-p},\eqno(2.7)$$
where $S_{i}(K,\cdot)$ is a regular Boel measure on $S^{n-1}$. The
measure $S^{n-1}(K,\cdot)$ is independent of the body $K$, and is
just ordinary Lebesgue measure, $S$, on $S^{n-1}$.
$S_{i}(B,\cdot)$ denotes the $i$-th surface area measure of the
unit ball in ${\Bbb R}^{n}$. In fact, $S_{i}(B,\cdot)=S$ for all
$i$. The surface area measure $S_{0}(K,\cdot)$ just is
$S(K,\cdot)$. When $i=0$, $S_{p,i}(K,\cdot)$ is written as
$S_{p}(K,\cdot)$ (see [26] and [27]). A fundamental inequality for
mixed $p$-quermassintegrals stats that: For $K,L\in{\cal
K}_{o}^{n}, p>1$ and $0\leq i<n-1$,
$$W_{p,i}(K,L)^{n-i}\geq W_{i}(K)^{n-i-p}W_{i}(L)^{p},\eqno(2.8)$$
with equality if and only if $K$ and $L$ are homothetic.
$L_{p}$-Brunn-Minkowski inequality for quermassintegrals
established by Lutwak [25]. If $K,L\in {\cal K}_{o}^{n}$ and $p\geq 1$ and $0\leq i\leq n$, then
$$W_{i}(K+_{p}L)^{p/(n-i)}\geq W_{i}(K)^{p/(n-i)}+W_{i}(L)^{p/(n-i)},\eqno(2.9)$$
with equality if and only if $K$ and $L$ are dilates or $L=\{o\}$.
Obviously, putting $i=0$ in (2.6), the mixed $p$-quermassintegrals
$W_{p,i}(K,L)$ become the well-known $L_{p}$-mixed volume
$V_{p}(K,L)$, defined by (see e.g. [26])
$$V_{p}(K,L)=\frac{1}{n}\int_{S^{n-1}}h_{L}^{p}dS_{p}(K,u).\eqno(2.10)$$

\vskip 10pt \noindent{\large \bf 3~ Orlicz multiple
mixed volumes}\vskip 10pt

In [18], the Orlicz multiple mixed volume was introduced as follows:

{\bf Definition 3.1}~ For $\varphi\in\Phi$, the Orlicz multiple mixed volume of convex bodies $K_{1},\cdots,K_{n},L_{n}$, denoted by $V_{\varphi}(K_{1},\cdots,K_{n},L_{n})$, defined by
$$V_{\varphi}(K_{1},\cdots,K_{n},L_{n})=\frac{1}{n}\int_{S^{n-1}}
\varphi\left(\frac{h_{L_{n}}}{h_{K_{n}}}\right)h_{K_{n}}dS(K_{1},\ldots,K_{n-1};u).\eqno(3.1)$$

When $K_{1}=\cdots=K_{n-i-1}=K$, $K_{n-i}=\cdots=K_{n-1}=B$,
$K_{n}=L$ and $L_{n}=K$, $V_{\varphi}(K_{1},\cdots,K_{n},L_{n})$
becomes the $i$-th Orlicz mixed volume $W_{\varphi,i}(K,L)$. When $L_{n}=K_{n}$
$V_{\varphi}(K_{1},\cdots,K_{n},L_{n})$ becomes the classical mixed volume
$V(K_{1},\cdots,K_{n}).$  When $K_{1}=\cdots=K_{n-1}=K$, $K_{n}=L$ and $L_{n}=K$, $V_{\varphi}(K_{1},\cdots,K_{n},L_{n})$
becomes the well-known Orlicz mixed volume $V_{\varphi}(K,L).$

{\bf The Orlicz-Aleksandrov-Fenchel inequality}~ {\it If $K_{1},\cdots,K_{n},L_{n}\in{\cal K}_{o}^{n}$, $1\leq r\leq n$ and ${\varphi}\in\Phi$, then} (see [18])
$$V_{\varphi}(K_{1},\cdots,
K_{n},L_{n})\geq V(K_{1},\cdots, K_{n-1},L_{n})\cdot
\varphi\left(\frac{\prod_{i=1}^{r}
V(K_{i}\ldots,K_{i},K_{r+1},\ldots,K_{n})^{\frac{1}{r}}}{V(K_{1},\cdots,
K_{n-1},L_{n})}\right).\eqno(3.2)$$

When $r=n$, $K_{1}=\cdots=K_{n-1}=K$, $K_{n}=L$ and
$L_{n}=K$, (3.2) becomes the following inequality which was
established by Gardner, Hug and Weil [28]. If $K$ is a convex body
containing the origin in its interior, $L$ is a convex body
containing the origin, and ${\varphi}\in\Phi$, then
$$V_{\varphi}(K,L)\geq V(K)\cdot
\varphi\left(\left(\frac{V(L)}{V(K)}\right)^{1/n}\right).\eqno(3.3)$$
If $\varphi$ is strictly convex, equality holds if and only if $K$
and $L$ are homothetic. Obviously, (3.2) in
special case yields the following result. If $K$ is a convex body
containing the origin in its interior, $L$ is a convex body
containing the origin, $0\leq i<n$ and ${\varphi}\in\Phi$, then
$$W_{\varphi,i}(K,L)\geq W_{i}(K)\cdot
\varphi\left(\left(\frac{W_{i}(L)}{W_{i}(K)}\right)^{1/(n-i)}\right).\eqno(3.4)$$
If $\varphi$ is strictly convex, equality holds if and only if $K$
and $L$ are homothetic, where $W_{i}(K)$ is the usual
quermassintegral of convex body $K$. (3.4) in special cases yield the following two well-known $L_{p}$-Minkowski type inequalities. $L_{p}$-Minkowski inequality for mixed
$p$-quermassintegrals established by Lutwak [26].
$$W_{p,i}(K,L)^{n-i}\geq W_{i}(K)^{n-i-p}W_{i}(L)^{p},\eqno(3.3)$$
for $p>1$ and $0\leq i\leq n$, with equality if and only if $K$
and $L$ are dilates or $L=\{o\}$. $L_{p}$-Minkowski inequality established by Firey
[29] . For $p>1$,
$$V_{p}(K,L)\geq V(K)^{(n-p)/n}V(L)^{p/n},\eqno(3.4)$$
with equality if and only if $K$ and $L$ are dilates or $L=\{o\}$.

\vskip 10pt \noindent{\large \bf 4~ Orlicz log-Aleksandrov-Fenchel inequality}\vskip 10pt

In the section, in order to prove the Orlicz-Aleksandrov-Fenchel inequality, we need to define some new mixed volume measures.

Associated with the convex bodies $L_{1},\cdots,L_{n-1}$ in ${\Bbb
R}^{n}$ is a unique positive Borel measure on $S^{n-1}$,
$S(L_{1},\ldots,L_{n-1};\cdot)$, call it the mixed area measure of
$L_{1},\ldots,L_{n-1}$, with the property that any convex body
$L_{n}$ one has the integral representation.
$$V(L_{1},\ldots,L_{n})=\frac{1}{n}\int_{S^{n-1}}h_{L_{n}}dS(L_{1},\ldots,L_{n-1};u).\eqno(4.1)$$

From (4.1), we introduce the mixed volume measure of convex bodies $L_{1},\ldots,L_{n}$.

{\bf Definition 4.1}~ (mixed volume measure) For $L_{1},\cdots,L_{n}\in{\cal K}^{n}_{o}$, the mixed volume measure of $L_{1},\ldots,L_{n}$, denoted by $dv(L_{1},\ldots,L_{n})$, defined by
$$dv(L_{1},\ldots,L_{n})=\frac{1}{n}h_{L_{n}}dS(L_{1},\cdots,L_{n-1};\cdot).\eqno(4.2)$$

From Definition 4.1, it is not difficult to find the following mixed volume probability measure.
$$d\bar{V}(L_{1},\ldots,L_{n})=\frac{1}{V(L_{1},\ldots,L_{n})}dv(L_{1},\ldots,L_{n}).\eqno(4.3)$$

For $p\geq 1$, $L_{p}$-multiple mixed volume of $L_{1},\cdots,L_{n},K_{n}$, denoted by $V_{p}(L_{1},\cdots,L_{n-1},K_{n},L_{n})$, defined by
$$V_{p}(L_{1},\cdots,L_{n-1},K_{n},L_{n})=\frac{1}{n}\int_{S^{n-1}}
\left(\frac{h_{K_{n}}}{h_{L_{n}}}\right)^{p}h_{L_{n}}dS(L_{1},\ldots,L_{n-1};u).\eqno(4.4)$$

From (4.4), we introduce $L_{p}$-multiple mixed volume measure of convex bodies $L_{1},\cdots,L_{n},K_{n}$ as follows.

{\bf Definition 4.2}~ (Orlicz multiple mixed volume measure) For $L_{1},\cdots,L_{n},K_{n}\in{\cal K}^{n}_{o}$, the Orlicz mixed volume measure of $L_{1},\ldots,L_{n},K_{n}$, denoted by $dv_{\varphi}(L_{1},\cdots,L_{n-1},K_{n},L_{n})$, defined by
$$dv_{\varphi}(L_{1},\cdots,L_{n-1},K_{n},L_{n})=\frac{1}{n}\cdot\varphi\left(\frac{h_{K_{n}}}{h_{L_{n}}}\right)\cdot h_{L_{n}}dS(L_{1},\ldots,L_{n-1};\cdot).\eqno(4.5)$$

From Definition 4.2, $L_{p}$-multiple mixed volume probability measure is defined by
$$d\bar{V}_{\varphi}(L_{1},\cdots,L_{n-1},K_{n},L_{n})=\frac{1}{V_{\varphi}(L_{1},
\cdots,L_{n-1},K_{n},L_{n})}dv_{\varphi}(L_{1},\cdots,L_{n-1},K_{n},L_{n}).\eqno(4.6)$$

Obviously, the mixed volume measure $dv(L_{1},\ldots,L_{n})$ is special case of the Orlicz multiple mixed volume measure. When $K_{n}=L_{n}$, we have
$$dv_{\varphi}(L_{1},\cdots,L_{n},L_{n})=dv(L_{1},\cdots,L_{n}).\eqno(4.7)$$

{\bf Theorem 4.3}~ (Orlicz log-Aleksandrov-Fenchel inequality) {\it If $L_{1},\ldots,L_{n},K_{n}\in {\cal K}^{n}_{o}$, $1\leq r\leq n$ and $\varphi\in \Phi$, then}
$$\int_{S^{n-1}}\ln\left(\varphi\left(\frac{h_{K_{n}}}{{h_{L_{n}}}}\right)\right)
d\bar{V}_{\varphi}(L_{1},\ldots,L_{n-1},K_{n},L_{n})\geq\ln\left(\varphi\left(\frac{\prod_{i=1}^{r}V(L_{i},\ldots,L_{i},L_{r+1},\ldots,L_{n-1},K_{n})^{1/r}}
{V(L_{1},\ldots,L_{n})}\right)\right).\eqno(4.9)$$

{\it Proof}~ From (4.2), (4.5) and (4.6), we have
$$\int_{S^{n-1}}\varphi\left(\frac{h_{K_{n}}}{h_{L_{n}}}\right)\ln\left(\frac{h_{K_{n}}}{h_{L_{n}}}\right)
dv(L_{1},\ldots,L_{n})=\int_{S^{n-1}}\ln\left(\frac{h_{K_{n}}}{h_{L_{n}}}\right)
dv_{\varphi}(L_{1},\ldots,L_{n-1},K_{n},L_{n}).\eqno(4.10)$$
Noting that
$$V_{\varphi}(L_{1},\ldots,L_{n-1},K_{n},L_{n})=\frac{1}{n}\int_{S^{n-1}}
\varphi\left(\frac{h_{K_{n}}}{h_{L_{n}}}\right)h_{L_{n}}dS(L_{1},\ldots,L_{n-1};u),$$
and from Lebesgue’s dominated convergence theorem, we obtain
$$\int_{S^{n-1}}\varphi\left(\frac{h_{K_{n}}}{h_{L_{n}}}\right)^{\frac{q}{q+n}}dv(L_{1},\ldots,L_{n})\rightarrow V_{\varphi}(L_{1},\ldots,L_{n-1},K_{n},L_{n})$$
as $q\rightarrow\infty,$ and
$$\int_{S^{n-1}}\varphi\left(\frac{h_{K_{n}}}{h_{L_{n}}}\right)^{\frac{q}{q+n}}\ln
\left(\frac{h_{K_{n}}}{h_{L_{n}}}\right)
dv(L_{1},\ldots,L_{n})\rightarrow\int_{S^{n-1}}\ln\left(\frac{h_{K_{n}}}{h_{L_{n}}}\right)
dv_{\varphi}(L_{1},\ldots,L_{n-1},K_{n},L_{n})$$
as $q\rightarrow\infty.$

Considering the function $f_{L_{1},\ldots,L_{n-1},K_{n},L_{n}}:[1,\infty]\rightarrow {\Bbb R}$, defined by
$$f_{L_{1},\ldots,L_{n-1},K_{n},L_{n}}(q)=\frac{1}{V_{\varphi}(L_{1},\ldots,L_{n-1},K_{n},L_{n})}\int_{S^{n-1}}
\varphi\left(\frac{h_{K_{n}}}{h_{L_{n}}}\right)^{\frac{q}{q+n}}dv(L_{1},\ldots,L_{n}).\eqno(4.11)$$

By calculating the derivative and limit of this function, we have
$$\frac{df_{L_{1},\ldots,L_{n-1},K_{n},L_{n}}(q)}{dq}=\frac{n}{(q+n)^{2}}\cdot\frac{1}{V_{\varphi}(L_{1},
\ldots,L_{n-1},K_{n},L_{n})}~~~~~~~~~~~$$
$$~~~~~~~~~~~~~~~~~~~~~~~~~~~\times\int_{S^{n-1}}\varphi\left(\frac{h_{K_{n}}}{h_{L_{n}}}\right)
^{\frac{q}{q+n}}\ln\left(\varphi\left(\frac{h_{K_{n}}}{h_{L_{n}}}\right)\right)
dv(L_{1},\ldots,L_{n}).\eqno(4.12)$$
and
$$\lim_{q\rightarrow\infty}f_{L_{1},\ldots,L_{n-1},K_{n},L_{n}}(q)=1.\eqno(4.13)$$

From (4.11), (4.12) and (4.13), and by using L'H\^{o}pital's rule, we have
\begin{eqnarray*}\lim_{q\rightarrow\infty}\ln\left(f_{L_{1},\ldots,L_{n-1},K_{n},L_{n}}(q)\right)^{q+n}
&=&-(q+n)^{2}\lim_{q\rightarrow\infty}\frac{\frac{\displaystyle df_{L_{1},\ldots,L_{n-1},K_{n},L_{n}}(q)}{\displaystyle dq}}
{f_{L_{1},\ldots,L_{n-1},K_{n},L_{n}}(q)}\\
&=&-\frac{n}{V_{\varphi}(L_{1},\ldots,L_{n-1},K_{n},L_{n})}\\
&\times&\lim_{q\rightarrow\infty}\frac{\int_{S^{n-1}}\varphi\left(\frac{h_{K_{n}}}{h_{L_{n}}}\right)
^{\frac{q}{q+n}}\ln\left(\varphi\left(\frac{h_{K_{n}}}{h_{L_{n}}}\right)\right)
dv(L_{1},\ldots,L_{n})}{f_{L_{1},\ldots,L_{n-1},K_{n},L_{n}}(q)}\\
&=&-\frac{n}{V_{\varphi}(L_{1},\ldots,L_{n-1},K_{n},L_{n})}\\
&\times & \int_{S^{n-1}}\varphi\left(\frac{h_{K_{n}}}{h_{L_{n}}}\right)
\ln\left(\varphi\left(\frac{h_{K_{n}}}{h_{L_{n}}}\right)\right)
dv(L_{1},\ldots,L_{n}).\end{eqnarray*}
Hence
$$\exp\left(-\frac{n}{V_{\varphi}(L_{1},\ldots,L_{n-1},K_{n},L_{n})}
\int_{S^{n-1}}\varphi\left(\frac{h_{K_{n}}}{h_{L_{n}}}\right)
\ln\left(\varphi\left(\frac{h_{K_{n}}}{h_{L_{n}}}\right)\right)
dv(L_{1},\ldots,L_{n})\right)$$
$$=\lim_{q\rightarrow\infty}(f_{L_{1},\ldots,L_{n-1},K_{n},L_{n}})^{q+n}~~~~~~~~~~~~~~~~~~~~~~~~~~~~~~~~~~~~~~~~~~~~$$
$$~~~~~~~~~~~~~~~~~~~~~=\lim_{q\rightarrow\infty}\left(\frac{1}
{V_{\varphi}(L_{1},\ldots,L_{n-1},K_{n},L_{n})}\int_{S^{n-1}}
\varphi\left(\frac{h_{K_{n}}}{h_{L_{n}}}\right)^{\frac{q}{q+n}}dv(L_{1},\ldots,L_{n})\right)^{q+n}.\eqno(4.14)$$
On the other hand, from H\"{o}lder's inequality
$$\left(\int_{S^{n-1}}\varphi\left(\frac{h_{K_{n}}}{h_{L_{n}}}\right)^{\frac{q}{q+n}}dv(L_{1},\ldots,L_{n})
\right)^{(q+n)/q}\left(\int_{S^{n-1}}dv(L_{1},\ldots,L_{n})\right)^{-n/q}$$
$$\leq\int_{S^{n-1}}\varphi\left(\frac{h_{K_{n}}}{h_{L_{n}}}\right)dv(L_{1},\ldots,L_{n})~~~~~~~~~~~~~~~~~~~~~~~~~~~~~$$
$$=V_{\varphi}(L_{1},\ldots,L_{n-1},K_{n},L_{n}).~~~~~~~~~~~~~~~~~~~~~~~~~~~~~~~~~~~\eqno(4.15)$$
From the equality condition of H\"{o}lder's inequality, it follows the equality in (4.15) holds if and only if
$h_{K_{n}}$ and $h_{L_{n}}$ are proportional. This yiels equality in (4.15) holds if and only if $K_{n}$ are $L_{n}$ homothetic.

Namely
$$\left(\frac{1}{V_{\varphi}(L_{1},\ldots,L_{n-1},K_{n},L_{n})}\int_{S^{n-1}}\varphi\left(\frac{h_{K_{n}}}{h_{L_{n}}}
\right)^{\frac{q}{q+n}}dv(L_{1},\ldots,L_{n})
\right)^{q+n}\leq\left(\frac{V(L_{1},\ldots,L_{n})}{V_{\varphi}(L_{1},\ldots,L_{n-1},K_{n},L_{n})}\right)^{n},$$
with equality if and only if $K_{n}$ are $L_{n}$ homothetic.

Hence
$$\exp\left(-\frac{n}{V_{\varphi}(L_{1},\ldots,L_{n-1},K_{n},L_{n})}
\int_{S^{n-1}}\varphi\left(\frac{h_{K_{n}}}{h_{L_{n}}}\right)
\ln\left(\varphi\left(\frac{h_{K_{n}}}{h_{L_{n}}}\right)\right)
dv(L_{1},\ldots,L_{n})\right)$$
$$\leq\left(\frac{V(L_{1},\ldots,L_{n})}{V_{\varphi}(L_{1},\ldots,L_{n-1},K_{n},L_{n})}\right)^{n},$$
with equality if and only if $K_{n}$ are $L_{n}$ homothetic.

That is
$$\frac{1}{V_{\varphi}(L_{1},\ldots,L_{n-1},K_{n},L_{n})}
\int_{S^{n-1}}\varphi\left(\frac{h_{K_{n}}}{h_{L_{n}}}\right)
\ln\left(\varphi\left(\frac{h_{K_{n}}}{h_{L_{n}}}\right)\right)
dv(L_{1},\ldots,L_{n})$$
$$\geq\ln\left(\frac{V_{\varphi}(L_{1},\ldots,L_{n-1},K_{n},L_{n})}
{V(L_{1},\ldots,L_{n})}\right),$$
with equality if and only if $K_{n}$ are $L_{n}$ homothetic.

Therefore
$$\int_{S^{n-1}}\ln\left(\varphi\left(\frac{h_{K_{n}}}{{h_{L_{n}}}}\right)\right)
d\bar{V}_{\varphi}(L_{1},\ldots,L_{n-1},K_{n},L_{n})\geq\ln\left(\frac{V_{\varphi}(L_{1},\ldots,L_{n-1},K_{n},L_{n})}
{V(L_{1},\ldots,L_{n})}\right),\eqno(4.16)$$
with equality if and only if $K_{n}$ are $L_{n}$ homothetic.

Further, on the right side of (4.16), by using the Orlicz-Aleksandrov-Fenchel inequality (3.2), (4.16) becomes (4.9).

This completes the proof. \hfill $\Box$

Unfortunately, precise equality for logarithmic Aleksandrov-Fenchel inequality (the second inequality in (4.9)) are also unknown in general, because the precise equality for the classical Aleksandrov-Fenchel inequality is unknown in general, see [19, p.376] for a full discussion.

{\bf Corollary 4.4}~ {\it If $K,L\in{\cal K}^{n}_{0}$ and $\varphi\in\Phi$, then
$$\int_{S^{n-1}}\ln\left(\frac{h_{K}}{h_{L}}\right)d\overline{W}_{\varphi,i}(L,K)
\geq\frac{1}{n-i}\ln\left(\frac{W_{i}(K)}{W_{i}(L)}\right).
\eqno(4.17)$$
each equality holds if and only if $K$ and $L$ are homothetic. Here $$dw_{\varphi,i}(L,K)=\frac{1}{n}\varphi\left(\frac{h_{K}}{h_{L}}\right)h_{L}dS_{i}(L,u),\eqno(4.18)$$
and $$d\overline{W}_{\varphi,i}(L,K)=\frac{1}{W_{\varphi,i}(L,K)}dw_{\varphi,i}(L,K),\eqno(4.19)$$
denotes its normalization.}

{\it Proof}~ This follows immediately from Theorem 4.3. \hfill $\Box$

\end{document}